\documentclass[11pt,  reqno]{amsart}%
\usepackage[french]{babel}
\usepackage[latin1]{inputenc}
\usepackage{graphics}
\usepackage{ulem}
\usepackage{hhline}
\usepackage{dsfont}
\usepackage{mathrsfs}
\usepackage{fancyhdr}
\usepackage{amsmath,  amssymb}
\usepackage{rotating}
\usepackage[latin1]{inputenc}
\usepackage[T1]{fontenc}
\usepackage{fancybox}
\usepackage{color}
\usepackage{colortbl}
\usepackage{setspace}
\usepackage{enumerate}
\usepackage[nice]{nicefrac}
\usepackage{amsthm}
\usepackage{multicol}
\usepackage{pifont}
\usepackage[french]{minitoc}
\usepackage{latexsym,  amsfonts}
\usepackage{amsthm}
\usepackage{varioref}
\usepackage{textcomp}
\usepackage{lmodern}
\usepackage{mathpazo}
\usepackage{euscript}
\usepackage[pdftex]{hyperref}
\numberwithin{equation}{section}
\makeatletter\@addtoreset{equation}{section}

\DeclareMathSymbol{\subsetneqq}{\mathbin}{AMSb}{36}
\newtheorem {theorem}{Théorème}[section]
\newtheorem {theo}{Theorem}[section]

\newtheorem {lemma}[theorem]{Lemme}

\newtheorem {remark}[theorem]{Remarque}
\newtheorem {remk}[theorem]{Remark}

\newcommand{\C}{\mathbb C}

\newcommand{\Pn}{\mathbb{P}^n(\mathbb C)}
\newcommand{\Ph}{\mathbb{P}^n(\mathbb H)}

\begin{document}
%================================================================================================
\author[]{A. Hafoud}
\title[]{Formules explicites  du noyau de la chaleur sur l'espace projectif quaternionique }
\maketitle

\begin{abstract}
 Dans cette Note, on donne une représentation intégrale  et le dévellopement assymptotique explicites  du noyau de la chaleur $ H_n(t;x,y)$ associé au Laplacien de Fubini-Study sur l'espace projectif quaternionique $\Ph$, en terme de la fonction théta  de Jacobi $\vartheta_2 $ et  des polynômes de Jacobi du type $ P_l^{(2n-1,1)}(\cos(2d))$.\\

\noindent {Abstract.} In this note we give an explicit integral representation  and  an expanssion for the heat kernel $ H_n(t;x,y)$ associated to  Fubini-Study Laplacians on quaternionic projective spaces $\Ph, n \geq 1$. This was possible by establishing a real integral representation formula for Jacobi polynomials of type  $P_l^{(2n-1,1)}(\cos(2d)) $.
\end{abstract}

\section{Abridged english version}

Let $ M=G/H$ be a rank one compact  riemannian symmetric space
(i.e., spheres, complex projective  spaces, quaternionic projective spaces),
$\Delta _M$ the Laplace-Beltrami operator on $M$ and $E_M(t;x,y)$ the
associated heat kernel solving the following heat equation on $M$:
\begin{align*} \left\{
     \begin{array}{ll}
        \frac{\partial}{\partial t}E_M(t;x,y)=\Delta_{M}E_M(t;x,y) ,  & \qquad t>0, \; x,y \in M \\
       \lim_{t\rightarrow 0^+}\int_{M}E_M(t;x,y)f(y)dy= f(x),         & \qquad f\in \mathcal{C}^\infty(M) \\
     \end{array}
   \right.. \end{align*}
Then  it is well known that the above heat kernel $E_M(t;x,y)$ depends only on $ t > 0 $ and on the geodesic distance $ d:=d(x,y)$ of the rank one compact symmetric space $M$, i.e.,  $E_M(t;x,y)= F(t,d)$. Hence it becomes natural to seek an explicit formula for the function $ F(t,d)$.For the case $ M= \Pn$, we have established the following explicit formulas for the heat kernel $Q_n(t,d):= Q_n(t;x,y)$ associated to  Fubini-Study Laplacian on complex projective space  $\Pn$ (see \cite{HafInt}):
\begin{align}
 Q_n(t;d) & = \frac{e^{n^2t}}{2^{n-2}{\pi}^{n+1}}
            \int_{d}^{\pi/2}  \frac{-d(cos(u))}{\sqrt{\cos^2(d)-\cos^2(u)}}
             \left[ -\frac 1{\sin(u)} \frac d{du}\right]^{n} \left(\theta_{n+1}(t;u)\right) \label{1}\\
Q_n(t;d) &=  \frac 1 {\pi^{n}} \sum_{l=0}^{+\infty} (2l+n)\frac{(l+n-1)!}{l!}
           e^{-4l(l+n)t} P_l^{(n-1,0)}(\cos(2d)) ,\label{2}
\end{align}
where $d:=d(x,y)$ is the geodesic distance on $\Pn$ and $\theta_{n+1}(t;u)$ is the function given by:
\begin{align*}
\theta_{n+1}(t;u) := \sum_{l=0}^{+\infty} e^{-4t(l+n/2)^2}\cos(2l+n)u; \quad  0 \leq d < \pi/2.
\end{align*}

In this Note, we give an integral representation and an expansion of the heat kernel
$H_n(t,d)$ associated to the Fubini-Study Laplacian on the quaternionic projective
space $\Ph$. Namely, we have:
 \begin{theo} Let $n \geq 1$. Then, the heat kernel $ H_n(t,d)$ associated to the
 Fubini-Study Laplacian on the quaternionic projective space $\Ph$ is given by the two following formulas:
\begin{align}
i) & \quad H_n(t;d)=  \frac{e^{(2n+1)^2t}}{2^{2n-2} \pi^{2n+1} }\int_d^{\pi/2}\frac{\sqrt{\cos^2(d)-\cos^2(u)}}{\cos^2(d)}\\
      &\qquad \qquad \qquad \qquad \qquad \qquad \times \left[ -\frac 1{\sin(u)} \frac d{du}\right]^{2n+1} \left[\theta_{2n+2}(t;u)\right]sin(u)du \label{3}\\
ii) & \quad H_n(t;d)=  \frac 1 {\pi^{2n}} \sum_{l=0}^{+\infty} (2l+2n+1)\frac{(l+2n)!}{(l+1)!}
           e^{-4l(l+2n+1)t} P_l^{(2n-1,1)}(\cos(2d)),   \label{4}
\end{align}
\end{theo}
where the function $\theta_{2n+2}(t;u)$  is given by:
\begin{align*}
\theta_{2n+2}(t;u) := \sum_{l=0}^{+\infty} e^{-4t(l+n+1/2))^2}\cos(2l+2n+1)u.
\end{align*}

\noindent {\bf Method of the proof:}
\begin{itemize}
 \item[(i)] relies essentially on an explicit integral representation of the heat kernel on the complex projective space $ \mathbb{P}^{2n+1}(\C)$  (see \cite{HafInt}), and
 \item[(ii)] is based on i) of Theorem 1.1 and  an real integral representation of the Jacobi polynomials of  type  $P_l^{(2n-1,1)}(x) $ .
 \end{itemize}
 To end this Abridged English version we make the following remarks:

 \begin{remk}
 The function $\theta_{2n+2}(t;x)$ above is related to the classical Jacobi function $\vartheta_2(z,\tau)$ (\cite[p.371]{Magnus}). More precisely we have:
  \begin{align*}
   \theta_{2n+2}(t;x)= \sum_{l=n}^{\infty} e^{-4t(l+1/2))^2}\cos(2l+1)x=\frac{1}{2} \vartheta_2\left(\frac{x}{\pi}, \frac{4it}{\pi}\right)-\sum_{l=0}^{n-1} e^{-4t(l+1/2))^2}\cos(2l+1)x.
  \end{align*}
 \end{remk}

 \begin{remk}
 Note that the kernels $ Q_n(t,d)$ and $ H_n(t,d)$ can be given by an unified formula (see French version).
 \end{remk}

\section{Introduction  et énnoncé des résultats }

Soit $\C$ le corps des nombres complexes et $ \mathds{P}^{2n+1}(\C)$ la variéte projective complexe  de $\C^{2n+2}$ et $ \mathds{H}=\C +\C j \approx \C^2$ le corps des nombres quaternioniques et $ \Ph$ la variéte projective quaternionique  de $\mathds{H}^{n+1}$, comme $ \C^{2n+2}$   peut s'identifier à $ \mathds{H}^{n+1}$, on peut alors considerer la projection $\pi$ donnée par:
\begin{align*}
\pi:  \mathbb{P}^{2n+1}(\C)  {\longrightarrow}\Ph . \end{align*}
\begin{align*}
 [z_0:z_1:......:z_{2n+1}]  {\longrightarrow} [z_0+z_1j:......:z_{2n}+z_{2n+1}j]  . \end{align*}
La varieté projective complexe $\mathds{P}^{2n+1}(\C)$ étant munie de sa métrique canonique de Fubini-Study $ds^2_{FS}$ et $\Ph$ est  munie de sa métrique canonique notée encore  $ ds^2_{FS}$ de sorte que $\pi$ soit  une submersion riemannienne dont toutes les fibres $(\mathbb{P}^{1}(\C);ds^2_{FS}) $ sont totalement géodésiques. Si   $\Delta_{\Pn} $  et $\Delta_{\Ph}$ désignent respectivement les Laplaciens sur $(\mathbb{P}^{2n+1}(\C);ds^2_{FS})$ et $(\Ph;ds^2_{FS})$, on sait que la relation d'entrelacement suivante: $ \pi ^* \circ \Delta_{\Ph}= \Delta_{\mathbb{P}^{2n+1}(\C)} \circ  \pi ^*  $ a  lieu (voir  \cite{BO}) pour la submersion riemannienne  donnée par la fibration de Hopf:
\begin{align}
(\mathbb{P}^1(\C);ds^2_{FS}) \hookrightarrow (\mathbb{P}^{2n+1}(\C);ds^2_{FS}) \stackrel{\pi}{\longrightarrow}
(\Ph;ds^2_{FS}).  \label{5}
\end{align}
Dans cette Note, on donne une représentation intégrale et un dévellopement assymptotique de ce qu'on appelle le noyau de diffusion sur $(\Ph, ds^2_{FS})$, i.e., $H_n(t;x,y) $ solution du problème de la chaleur associé à  $\Delta_{\Ph} $  sur l'espace $\Ph$:
\begin{align*}
\left\{
     \begin{array}{ll}
       \frac{\partial}{\partial t}H_n(t;x,y)= \Delta_{\Ph}H_n(t;x,y) ,  &\quad t>0, x,y \in \Ph \\
       \lim_{t\rightarrow 0^+}\int_{ \Ph}H_n(t;x,y)f(y)dy= f(x), x \in \Ph , &\quad f\in \mathcal{C}^\infty( \Ph ) \\
     \end{array}
     \right. .\qquad
      (H_n)
\end{align*}
Pour énoncer les résultats principaux de cette Note, on fixera quelques notations qui seront utilisées par la suite.
  La distance géodésique $d_{FS}(x,y)$  sur $\mathbb{P}^{2n+1}(\C)$ est donnée par:
\begin{align}
\cos(d_{{FS}}(x,y))= \frac{\mid \sum_{i=0}^{2n+1}\overline{x_i}y_i  \mid}{\mid x \mid  \mid y\mid},\quad  x, y \in  \mathbb{P}^{2n+1}(\C) .\label{7}
\end{align}
La distance $d_{FS}(x,y)$  sur $\Ph$ est définie  par:
\begin{align}
 \cos(d_{FS}(x,y))= \frac{\mid \sum_{i=0}^{n}\overline{x_i}y_i  \mid}{\mid x \mid  \mid y\mid},\quad  x, y  \in  \Ph  .\label{8}
 \end{align}
Dans la suite on notera tout simplement $ d(x,y)$. Pour  $ \alpha, \beta > -1/2$, $ l \in \mathds{N} $,  $ P_l^{(\alpha,\beta)}(x)$ désigne le polynôme de Jacobi de degrés $l$, $C_l^1(x)$ désigne le polynôme de Gegenbauer qui vérifie la relation $C_l^1(\cos(\theta))= \frac{\sin(l+1)\theta}{\sin(\theta)}$. On note par $\pi^*$ l'application  qui à $ f$  associe $ \pi^*(f)=f\circ \pi $. On désigne par $L $ l'opérateur  différentiel $ L= -\frac 1{\sin(u)} \frac d{du}$ et  $ L^m$ le composé de $L$ m fois et par $ \theta_{2n+2}(t,u)$, $ \Psi_{m}(t,u)$   les  fonctions données réspectivement par:
\begin{align}
\theta_{2n+2}(t,u) := \sum_{l=0}^{+\infty} e^{-4t(l+n+1/2))^2}\cos(2l+2n+1)u, \quad \Psi_{m}(t,u)=\sin(u)L^{m} [\theta_{2n+2}(t,u)]. \label{9}
\end{align}
Le résultat principal de cette Note s'énonce comme suit:

\begin{theorem} \label{thm1}
Soit $ n  \in \mathds{N}^* $. Ecrivons $H_n(t,d)=H_n(t;x,y) $ avec $d:=d(x,y)$ le noyau de la chaleur associé  au Laplacien de  Fubini-Study  sur l'espace projectif quaternionique  $ \Ph$, o$\grave{u}$ $d$ est la distance de Fubini-Study sur $\Ph$. Alors, pour tout $ t > 0 $, le  noyau  $ H_n(t,d)$ admet les deux représentations suivantes:
\begin{align}
i) & \quad H_n(t;d)=  \frac{e^{(2n+1)^2t}}{2^{2n-2} \pi^{2n+1} }\int_d^{\pi/2}
\frac{\sqrt{\cos^2(d)-\cos^2(u)}}{\cos^2(d)} \Psi_{2n+1}(t,u) du
\label{10}\\
ii) & \quad H_n(t;d)=  \frac 1 {\pi^{2n}} \sum_{l=0}^{+\infty} (2l+2n+1)\frac{(l+2n)!}{(l+1)!}
           e^{-4l(l+2n+1)t} P_l^{(2n-1,1)}(\cos(2d)), \label{11}
\end{align}
o$\grave{u}$ la fonction  $ \Psi_{2n+1}(t,u)$ est donnée par:
 \begin{align*}
  \Psi_{2n+1}(t,u):=  \sin(u)\left( -\frac 1{\sin(u)} \frac d{du} \right)^{2n+1} \big( \sum_{l=0}^{+\infty} e^{-4t(l+n+1/2))^2}\cos(2l+2n+1)u \big).
 \end{align*}
 \end{theorem}

Avant de donner une esquisse de la preuve du théorème ci-dessus, on mentionne les remarques suivantes:

\begin{remark}
Les deux  noyaux de la chaleur $ H_n(t;d)$ et $ Q_n(t;d)$ peuvent s'unifier dans une seule formule . En effet: Soient $ \mathds{F}=\C$ ou   $\mathds{H}$, et  soit $ 2k$  la dimension réelle de $\mathds{F}$ considéré comme espace vectoriel  sur $ \mathds{R}$  ($k= 1$ ou $2$) et $M_{n,k}=
 \mathds{P}^{n}(\mathds{F})$ et $\Delta_{n,k}$ le Laplacien de Beltrami sur la variété  $ M_{n,k}$ et $E_{n,k}(t,d) $ le  noyau de la chaleur associé à $\Delta_{n,k}$. Alors $E_{n,k}(t,d) $   se met sous l'une des deux formes:
  \begin{align}
 i)& \quad E_{n,k}(t;d)= \frac 1 {\pi^{kn}} \sum_{l=0}^{+\infty} (2l+k(n+1)-1)\frac{(l+k(n+1)-2)!}{(l+k-1)!}
           e^{-4l(l+k(n+1)-1)t} P_l^{(kn-1,k-1)}(\cos(2d)) \label{12} \\
ii)& \quad E_{n,k}(t;d)  = \frac{c(n,k) e^{(k(n+1)-1)^2t}}{ \cos(d)^{2(k-1)} } \int_{d}^{\pi/2}  \frac{-d(cos(u))} { (\cos^2(d)-\cos^2(u))^{3/2-k}}
 L^{k(n+1) -1} \left(\theta_{k(n+1)}(t;u)\right) \label{13}
 \end{align}
avec
 \begin{align*}
 c(n,k) &=\frac{1}{2^{kn-2}  {\pi}^{kn+1 } } ,\quad  \theta_{k(n+1)}(t;u) = \sum_{l=0}^{+\infty} e^{-4t(l+\frac{k(n+1)-1}{2})^2}\cos(2l+k(n+1) -1)u.
 \end{align*}

 \end{remark}

 \begin{remark}
Le  developpement du noyau de la chaleur $ H_n(t;d)$ sur $\Ph$ donné dans ii) du Théorème \ref{thm1} donne en effet la forme explicite du developpement établit, en terme de fonctions zonales, par (\cite[Remarque 2.1, P.270]{Benabdelah}) du noyau de la chaleur de l'espace projectif quaternionique $\Ph$ vu comme espace homogéne .
\end{remark}

\section{Méthode  de démonstration du Théorème \ref{thm1}}

Preuve de i). Pour commencer, rappelons tout d'abord que l'on a la fibration de Hopf:
\begin{align}
(\mathbb{P}^1(\C);ds^2_{FS}) \hookrightarrow (\mathbb{P}^{2n+1}(\C);ds^2_{FS}) \stackrel{\pi}{\longrightarrow}
(\Ph;ds^2_{FS}).
\label{14}
\end{align}
  Soit $ u(t,x) $  la solution de l'équation de la chaleur associée à  $\Delta_{\Ph}$  sur $\Ph$:
\begin{align}\left\{
     \begin{array}{ll}
         \frac{\partial}{\partial t}u(t,x) =\Delta_{\Ph}u(t,x) , & \quad t>0, x \in \Ph \\
       u(0,x)= f_0(x), & \quad f_0 \in \mathcal{C}^\infty( \Ph ) \\
     \end{array}
   \right.  \qquad \label{15}
\end{align}
Comme on a la formule d'entrelacement $ \pi ^* \circ \Delta_{\Ph}= \Delta_{\mathbb{P}^{2n+1}(\C)} \circ  \pi ^*  $ alors $ (\pi ^*u)(t,x)$ est la solution de l'équation de la chaleur associée au Laplacien $ \Delta_{\mathbb{P}^{2n+1}(\C)}$ sur l'éspace projectif complexe $\mathbb{P}^{2n+1}(\C)$, donc en tenant  compte de la  formule integrale  explicite du noyau de la chaleur  $Q_{2n+1}(t; d(x,y))$ sur $\mathbb{P}^{2n+1}(\C)$  donnée dans l'équation (1) (voir \cite{HafInt}), on peut écrire :
\begin{align}
(\pi ^*u)(t,x)= \int_{\mathbb{P}^{2n+1}(\C)} Q_{2n+1}(t; d(x,y))(\pi ^*f_0)(y)dy \label{15b}
 \end{align}
o$\grave{u}$  le noyau  $Q_{2n+1}(t;d(x,y))$ est donné par:
\begin{align}
Q_{2n+1}(t;d)  = \frac{e^{(2n+1)^2t}}{2^{2n-1} \pi^{2n+2} } \int_{d}^{\pi/2}
\frac{1}{\sqrt{\cos^2(d)-\cos^2(u)}}
\Psi_{2n+1}(t;u)du. \label{16}
\end{align}
En introduisant la fonction de Haiveside $Y$, le noyau $Q_{2n+1}(t;d)$  peut se mettre sous la forme:
\begin{align}
Q_{2n+1}(t;d)= \frac{e^{(2n+1)^2t}}{2^{2n-1} \pi^{2n+2}}\int_{0}^{\pi/2}
\frac{Y(u-d)}{\sqrt{\cos^2(d)-\cos^2(u)}}\Psi_{2n+1} (t;u)du. \label{17}
\end{align}
Alors, en  remplaçant  dans l'équation  (3.3) le noyau  $Q_{2n+1}(t;d)$  par son expression ci-dessus  et en permuttant les intégrales,  la solution $(\pi ^*u)(t,x) $  peut se mettre sous la forme:
\begin{align}
u(t,\pi(x)) = \frac{e^{(2n+1)^2t}}{2^{2n-1} \pi^{2n+2} } \int_0 ^{\pi/2} T(u,x)\Psi_{2n+1} (t,u) du \label{18}
\end{align}
o$\grave{u}$   $ T(u,x)$ est l'intégrale suivante:
\begin{align}
T(u,x)= \int_{\mathbb{P}^{2n+1}(\C)} \frac{Y(u-d(x,y))}{\sqrt{\cos^2(d(x,y))-\cos^2(u)}}f_0(\pi(y)) dy .\label{19}
\end{align}
Pour $ x=O=[1:0:....0] \in \mathbb{P}^{2n+1}(\C)$,  en utilisant les m\^emes téchniques utilisées dans (\cite{BO}),   on  ramene l'intégrale $T(u,O)$ à une intégrale sur $\Ph$. En effet: si on désigne par $ B_u(O)$ la boule  de  $\mathbb{P}^{2n+1}(\C)$ de centre $O=[1:0:0....0]$ et de rayon $ u$, l'intégrale $T(u,O)$ se met sous la forme:
\begin{align}
T(u,O)= \int_{y \in B_u(O)} \frac{f_0(\pi(y))}{\sqrt{\cos^2(d(O,y))-\cos^2(u)}} dy \label{20}
\end{align}
et en utilisant le théorème de Fubini on a:
\begin{align}
T(u,O)= \int_{z \in B_u(\pi(O))}  \left( \int_{y \in  \pi ^{-1}(z)} \frac{dvol(\pi
^{-1}(z))(y)}{\sqrt{\cos^2(d(O,y))-\cos^2(u)}} \right) f_0(z) dz  \label{21}
\end{align}
et pour réduire  l'intégrale $T(u,O)$ on  paramétrise un voisinage de $ O=[1:0:.....:0] \in \Ph$ par:
\begin{align}
z=[1:\tan(\phi).q ]; \quad \phi \in [0,\pi /2[ ; \quad q \in S^{4n-1} \subset \mathds{H}^n \label{22}
\end{align}
alors, tout  élément $y $  de  $ \pi ^{-1}(z)$ est de la forme :
\begin{align}
y=[\cos(v):\sin(v)e^{i\theta}: \tan(\phi)q(\cos(v)+j\sin(v)e^{i\theta})]; \quad  \theta \in [0,2\pi]; v \in [0,\pi
/2] \label{23}
\end{align}
et la fibre $\pi ^{-1}(z)$ s'identifie à $\mathbb{P}^{1}(\C)$ dont les variables  $(v,\theta)$  sont les   coordonnées
radiale et angulaire respectivement, et l'on a :
\begin{align}
dvol(\pi ^{-1}(z))(y)= \sin(2v)dv d\theta/2 \label{24}
\end{align}
et en utilisant les formules données par les équations \eqref{7} et \eqref{8}  on a:
\begin{align}
\cos^2(d(O,y))=\cos^2(\phi)\cos^2(v) ,\quad    \cos^2(d(O,z))=cos^2(\phi),  \label{25}
\end{align}
et par un calcul direct on obtient:
\begin{align}
\int_{y \in  \pi ^{-1}(z)} \frac{dvol(\pi ^{-1}(z))(y)  }{\sqrt{\cos^2(d)-\cos^2(u)}}=
     2\pi   \frac{\sqrt{\cos^2(\phi)-\cos^2(u)}}{\cos^2(\phi)} .\label{26}
     \end{align}
L'intégrale  $ T(u,O)$   se met alors sous la forme:
\begin{align}
T(u,O)=  2\pi  \int_{z \in \Ph} Y(u-d(O,z))
\frac{\sqrt{\cos^2(d(O,z))-\cos^2(u)}}{\cos^2(d(O,z))} f_0(z) dz .\label{27}
\end{align}
Par conséquent,  la solution $ u(t,x)$ du probléme $ (H_n)$  évaluée
en $ O=[1:0:0.....0] \in \mathds{P}^n(\mathds{H})$ s'écrit sous la forme:
 \begin{align}
 u(t,O) = \frac{e^{(2n+1)^2t}}{2^{2n-2} \pi^{2n+1} }  \int_{ \Ph  } \int_{d(O,z)} ^{\pi/2} \frac{\sqrt{\cos^2(d(O,z))-\cos^2(u)
 }}{\cos^2(d(O,z))} \Psi_{2n+1} (t,u) f_0(z) du dz .\label{28}
 \end{align}
Comme le noyau  $H_{n}(t,x,y)$  ne dépend  que de la distance géodèsique $ d:= d(x,y)$ on conclut que :
 \begin{align}
 H_{n}(t; d)= \frac{e^{(2n+1)^2t}}{2^{2n-2} \pi^{2n+1}}\int_{d} ^{\pi/2}
 \frac{\sqrt{\cos^2(d)-\cos^2(u)}}{\cos^2(d)} \Psi_{2n+1}(t,u)du, \label{29}
 \end{align}
D'ou la formule i) du théorème.\\

 Pour la preuve  de ii) on  remarque que:
\begin{align}
\Psi_{2n+1}(t,u)= \sum_{l=0}^{+\infty} e^{-4t(l+n+1/2)^2} (2l+2n+1) \sin(u) L^{2n} \left(
C_{2l+2n}^1(\cos(u))\right) .\label{30}
\end{align}
Le noyau  $H_{n}(t; d)$ se met alors  sous la forme:
\begin{align}
 H_{n}(t;d)  &= \frac{e^{(2n+1)^2t}}{2^{2n-2}\pi^{2n+1}}\sum_{l=0}^{+\infty} e^{-4t(l+n+1/2)^2} (2l+2n+1) \label{31} \\
 &   \times \int_{d}^{\pi/2} \frac{\sqrt{\cos^2(d)-\cos^2(u)}}{\cos^2(d)} L^{2n}
  \left( C_{2l+2n}^1(\cos(u)) \right) sin(u)du.  \nonumber
\end{align}

Pour donner la formule ii)  du théorème on a besoin du Lemme suivant:

\begin{lemma}
Soient $ n \in \mathbb{N}^*, l \in \mathbb{N}, 0 \leq d <\pi/2 $,  alors les polyn\^omes de Jacobi du type $ P_l^{(2n-1,1)}(cos(2d)))$ admettent une représentation intégrale réelle en terme des polyn\^omes de Gegenbauer  $ C_{2l+2n}^1(\cos(u))$ et l'on a la formule suivante:
 \begin{align*}
 \int_{d} ^{\pi/2} \frac{\sqrt{\cos^2(d)-\cos^2(u)}}{\cos^2(d)}& L^{2n} \left( C_{2l+2n}^1(\cos(u)) \right)
 sin(u)du \\ &= \frac{2^{2n-2}\pi (l+2n)!}{(l+1)!}P_{l}^{(2n-1,1)}(cos(2d)).
 \end{align*}
\end{lemma}
\noindent En  utilisant le lemme ci-dessus,  le noyau $H_{n}(t,d)$ se met sous la forme ii) du théorème:
\begin{align}
H_{n}(t;d) =\frac 1 {\pi^{2n}} \sum_{l=0}^{+\infty} (2l+2n+1)\frac{(l+2n)!}{(l+1)!}
           e^{-4l(l+2n+1)t} P_l^{(2n-1,1)}(\cos(2d)). \label{32}
 \end{align}
 Ceci termine l'esquisse de la preuve du Théorème \ref{thm1}.\\

\noindent {\bf Idée de la preuve du Lemme:}
Les polyn\^omes de Jacobi du type $ P_{l+1}^{(2n-2,0)}(2t^2-1)$ admettent une représentation intégrale réelle (voir \cite{HafInt}):
\begin{align}
P_{l+1}^{(2n-1,0)}(\cos(2d)) = \frac{2 .(l+1)!(2n-2)!}{\pi (l+2n-1)!}  \int_d^{\pi/2}
 \frac{-d(\cos(u)}{\sqrt{\cos^2(d)-\cos^2(u)}} C_{2l+2}^{2n-1}(cos(u)). \label{33}
 \end{align}
 En appliquant des changements de variables  dans la formule donnée dans l'equation \eqref{33}  on obtient:
\begin{align}
P^{(2n-2,0)}_{l+1}(2t^2-1) = \frac{2. (l+1)!(2n-2)!}{\pi(l+2n-1)!}
\int^1_0  \frac{1}{\sqrt{1 -u^2}} C^{2n-1}_{2l+2}(ut)du.  \label{34}
\end{align}
Par dérivation par rapport à $t$ et en faisant une intégration par partie, puis en faisant des changements de variables inverses, on établit la formule du Lemme.\\

{\bf Remerciements:} Je tiens à remercier Mr Ahmed Intissar de la Faculté des Sciences de Rabat-Agdal pour les fructueuses discussions qu'il m'a accordées.

\end{document}